\documentclass[11pt]{amsart}
\usepackage{amsmath}
\usepackage{mathrsfs}
\def\pair#1{\langle #1 \rangle}
\def\Exp#1{\exp \left[#1\right]}
\theoremstyle{plain}
    \newtheorem{theorem}{Theorem}[section]
    \newtheorem{proposition}[theorem]{Proposition}
    \newtheorem{lemma}[theorem]{Lemma}
    \newtheorem{corollary}[theorem]{Corollary}
    \newtheorem{remark}[theorem]{Remark}
\theoremstyle{definition}
    \newtheorem{definition}[theorem]{Definition}

\newcommand{\bbZ}{\mathbb{Z}}
\newcommand{\bbC}{\mathbb{C}}
\newcommand{\bbQ}{\mathbb{Q}}

\def\col{\mathrm{col}}
\def\ol{\overline}
\def\cO{\mathcal{O}}
\def\bsK{{\boldsymbol{K}}}
\def\frp{\mathfrak{p}}

\begin{document}
\noindent

%--- The Title ------------------------------------------------------------------------------------------------------------------------
\title[Kronecker limit formulas]{The Kronecker limit formulas via the distribution relation}
\author{Kenichi Bannai}
\address{Department of Mathematics, Keio University, 3-14-1 Hiyoshi, Kouhoku-ku, Yokohama 223-8522, Japan}
\email{bannai@math.keio.ac.jp}
\author{Shinichi Kobayashi}
\address{Graduate School of Mathematics, Nagoya University, 
                           Furo-cho Chikusa-ku, Nagoya  464-8602, Japan}
\email{shinichi@math.nagoya-u.ac.jp}
\footnote{The 2000 Mathematics Subject Classification: 11M35, 11M36, 11S80}
\begin{abstract}
	In this paper, we give a proof of the classical Kronecker limit formulas using the distribution relation of the Eisenstein-Kronecker series.
	Using a similar idea, we then prove $p$-adic analogues of the Kronecker limit formulas for the 
	$p$-adic Eisenstein-Kronecker functions defined in our previous paper.
\end{abstract}

\maketitle

%%%%%%%%%%%%%%%%%%%%%%%%%%%%%%%%%%%%%%%%%%%%%%%%%%%%%%%
%
%
%
\section{Introduction}
%
%
%
%%%%%%%%%%%%%%%%%%%%%%%%%%%%%%%%%%%%%%%%%%%%%%%%%%%%%%%

In these notes, we will give a proof of the classical first and second Kronecker limit formulas concerning the
limit of values of Eisenstein-Kronecker-Lerch series.  Our proof is based on the distribution relation of the
Eisenstein-Kronecker-Lerch series.
Using a similar idea, we then prove $p$-adic analogues of these formulas 
for the $p$-adic Eisenstein-Kronecker functions defined in our previous paper \cite{BFK}.

Let $\Gamma \subset \bbC$ be a lattice.   We define a pairing for $z$, $w \in \mathbb{C}$ by 
$\pair{z,w}_\Gamma :=  \Exp {(z \overline{w}  - w \overline{z})/A(\Gamma)}$, where $A(\Gamma)$ is the area of the fundamental 
domain of $\Gamma$ divided by $\boldsymbol\pi = 3.1415\cdots$.
Then for an integer $a$ and a fixed $z_0$, $w_0 \in \bbC$,
the Eisenstein-Kronecker-Lerch series is defined as
\begin{equation*} \label{equation: definition of Eisenstein-Kronecker*}
	K^*_a(z_0,w_0,s; \Gamma) = {\sum}^*_{\gamma \in \Gamma} 
    		\frac{(\overline{z}_0 + \overline{\gamma})^a}{|z_0 + \gamma|^{2s}} \pair{\gamma, w_0}_\Gamma,
\end{equation*}
where $\sum^*$ denotes the sum taken over all $\gamma\in\Gamma$
except for $\gamma=-z_0$ if $z_0 \in \Gamma$.  The above series converges for
$\operatorname{Re}(s) > a/2+1$, but one may give it meaning for general $s$ by analytic continuation.
In what follows, we omit the $\Gamma$ from the notations if there is no fear of confusion.

We let $\theta(z)$ be the reduced theta function on $\bbC/\Gamma$ associated to the divisor $[0] \subset \bbC/\Gamma$, normalized so that $\theta'(0) =1$.  
Then the Kronecker limit formulas are given as follows.

\begin{theorem}[Kronecker limit formulas]\label{thm: KLF}
	Let  $c$ be the Euler constant 
	$
		c:=\lim_{n\rightarrow \infty} \left(1+\frac{1}{2}+ \cdots +\frac{1}{n}-\log n\right),
	$ 
	and let $\Delta$ be the discriminant of  $\Gamma$ defined as $\Delta:=g_2^3-27g_3^2$, where 
	$
		g_k:=\sum_{\gamma \in \Gamma \setminus \{0\}} \gamma^{-2k}.
	$
 	Then we have the following.
	\begin{enumerate}
		\item The first limit formula
		$$
			\lim_{s \rightarrow 1}\left(A K^*_0(0, 0, s)-\frac{1}{s-1}\right)
			=-\frac{1}{12}\log|\Delta|^2-2\log A+2c.
		$$
		\item For  $z \notin \Gamma$, the second limit formula
		$$
			A K^*_0(0,z,1)=-\log |\theta(z)|^2+\frac{|z|^2}{A}-\frac{1}{12}\log|\Delta|^2.
		$$	
	\end{enumerate}
\end{theorem}

Numerous proofs exist for this classical formula, and many of the proofs rely on arguments concerning the 
moduli space.  We give a proof of the above theorem, valid for a fixed lattice $\Gamma \subset \bbC$,
using the Kronecker theta function.   
As in the original proof by Kronecker, we first prove the second limit formula, and then 
deduce the first limit formula from the second.
The key in our proof is the distribution relation for the Eisenstein-Kronecker function.

Our view of understanding the Kronecker limit formulas in terms of the Kronecker theta function
and the distribution relation allows us to prove the following $p$-adic analogues of Theorem \ref{thm: KLF}.
Suppose now that $\Gamma$ corresponds to a period lattice corresponding to the invariant differential $\omega = dx/y$
of an elliptic curve $E: y^2 = 4x^3 - g_2 x - g_3$ with complex multiplication by
the ring of integers $\cO_\bsK$ of an imaginary quadratic field $\bsK$.  We assume in addition that $E$ is defined over $\bsK$,
and that the model above has good reduction at the primes above $p \geq 5$.
We fix a branch of the $p$-adic logarithm, which is a homomorphism $\log_p : \bbC_p^\times \rightarrow \bbC_p$.
In the paper \cite{BFK}, we introduced the $p$-adic Eisenstein-Kronecker series $E^\col_{a,b}(z)$ as a Coleman function
for integers $a$, $b$ such that $b \geq 0$.  This function is a $p$-adic analogue of the Kronecker double series
$$
	K^*_{a+b}(0,z,b) = {\sum}^* _{\gamma\in\Gamma}\frac{\ol\gamma^{a+b}}{|\gamma|^{b}} \pair{\gamma,z}.
$$
In this paper, we let
$$
	K^\col_{a+b}(0,z,b) := E^\col_{a,b}(z)
$$
to highlight the analogy.  Then in analogy with Theorem \ref{thm: KLF} (ii), we have the following.

\begin{theorem}[$p$-adic Kronecker second limit formula]\label{thm: pKLF}
	For any prime $p \geq 5$ of good reduction, we have the second limit formula
	$$
		 K^\col_0(0,z,1)= - \log_p \theta(z) - \frac{1}{12}\log_p \Delta,
	$$
	where $\log_p \theta(z)$ is a certain $p$-adic analogue of the function $\log|\theta(z)| - |z|^2/A$ 
	defined in Definition \ref{def: log-p-theta} using the reduced theta function $\theta(z)$ and the branch of our $p$-adic logarithm. 
\end{theorem}

The $p$-adic analogues of Kronecker second limit formula were previously investigated by Katz \cite{Ka} and de Shalit \cite{dS}
in the context of $p$-adic $L$-functions when $p$ is a prime of good ordinary reduction. 
Our formulation via $p$-adic Eisenstein-Kronecker series gives a direct $p$-adic analogue, and is valid even for supersingular $p$.

When $p \geq 5$ is a prime of good ordinary reduction, we defined 
in \cite{BK1} \S 3.1 a two-variable $p$-adic measure $\mu:=\mu_{0,0}$ on $\bbZ_p \times \bbZ_p$ interpolating Eisenstein-Kronecker numbers,
or more precisely, the values $K^*_{a+b}(0,0,b)/A^a$ for $a$, $b \geq 0$.
We define the $p$-adic Eisenstein-Kronecker-Lerch series by
$$
	K^{(p)}_{a}(0, 0, s) := \int_{\bbZ^\times_p \times \bbZ^\times_p} \pair{x}^{s-1} \omega(y)^{a-1}\pair{y}^{a-s} d \mu(x,y)
$$
for any $s \in \bbZ_p$, where $\pair{-} : \bbZ_p^\times \rightarrow \bbC_p^\times$ is given as the composition
$\bbZ_p^\times \rightarrow 1 + p\bbZ_p \hookrightarrow \bbC_p^\times$ and
$\omega: \bbZ_p^\times \rightarrow \mu_{p-1}$ is the Teichm\"uller character.  
Then an argument similar to the proof of Theorem \ref{thm: KLF} (i) gives the following.

\begin{proposition}\label{pro: pKLF}
 	Suppose $p \geq 5$ is  a prime of good ordinary reduction.  Then
	$$
		\lim_{s\rightarrow 1} K^{(p)}_{0}(0,0,s) = \Omega_p^{-1}\left( 1 - \frac{1}{p} \right) \log_p \ol\pi
	$$
	where $\Omega_p$ is a $p$-adic period of the formal group of $E$. 
\end{proposition}

The proof of the above proposition is similar to that of the proof of Theorem \ref{thm: KLF}.  However, due to the
existence of a trivial zero for the function $K^{(p)}_0(0,0,s)$ at $s=1$, the analogy with the classical case is not
perfect.  See Remark \ref{rem: not pKLF} for details.

%\setcounter{tocdepth}{1}
%\tableofcontents

%%%%%%%%%%%%%%%%%%%%%%%%%%%%%%%%%%%%%%%%%%%%%%%%%%%%%%%
%
%
%
\section{Kronecker limit formulas}
%
%
%
%%%%%%%%%%%%%%%%%%%%%%%%%%%%%%%%%%%%%%%%%%%%%%%%%%%%%%%

%%%%%%%%%%%%%%%%%%%%%%%%%%%%%%%%%%%%%%%%%%%%%%%%%%%%%%%
%
\subsection{Kronecker's Theorem}
%
%%%%%%%%%%%%%%%%%%%%%%%%%%%%%%%%%%%%%%%%%%%%%%%%%%%%%%%

In this section, we recall the definition of Eisenstein-Kronecker-Lerch series and the Kronecker theta function.
Then we will state Kronecker's theorem giving the relation between the two.  All of the results are contained in 
\cite{We1}.   Se also \cite{BK1} or \cite{BKT}.

We fix a lattice $\Gamma$ in $\mathbb{C}$ and let $A$ be the area of the fundamental domain of $\Gamma$ divided by $\pi$.  
Let $a$ be an integer $\geq 0$.  
For a fixed $z_0$, $w_0 \in \bbC$,  we let
$\theta^*_a(t,z_0,w_0)$ be the function
$$
	\theta^*_a(t,z_0,w_0) = {\sum_{\gamma \in \Gamma}}^* \exp(- t|z_0 + \gamma|^2/A) \pair{\gamma, w_0} (\overline z_0 +\overline \gamma)^a,
$$
where $\sum^*$  means the sum taken over all $\gamma \in \Gamma$ other than $-z_0$ if $z_0$ is in $\Gamma$.
Furthermore, we let
$$
	I_a(z_0, w_0, s) := \int_{1}^\infty\theta^*_a(t,z_0,w_0) t^{s-1} dt. 
$$
Then we have 
\begin{multline}\label{eq: integral expression}
	A^s \Gamma(s) K^*_a(z_0, w_0, s) = I_a(z_0, w_0, s)
	 - \frac{\delta_{z_0, a}}{s} \pair{w_0,z_0}\\ 
	+ I_{a} (w_0, z_0, a+1-s) \pair{w_0,z_0}  + \frac{\delta_{w_0, a}}{s-1},
\end{multline}
where $\delta_{a,x} = 1$ if $a =0$ and $x \in \Gamma$, and $\delta_{a,x} = 0$ otherwise.   The above integral
expression gives the meromorphic continuation of $K^*_a(z_0, w_0, s)$ to the whole complex plane,
and also the functional equation. 
    
We next review the definition of the Kronecker theta  function.  We let $\theta(z)$ be the reduced theta function
on $\bbC/\Gamma$ associated to the divisor $[0] \in \bbC/\Gamma$, normalized so that $\theta'(0) =1$.  This function 
may be given explicitly in terms of the Weierstrass $\sigma$-function
\begin{equation*}\label{product expansion of sigma}
	\sigma(z) : = z \prod_{\gamma \in \Gamma \setminus \{ 0 \}}
	\left( 1 - \frac{z}{\gamma} \right) \Exp{ \frac{z}{\gamma} + \frac{z^2}{2 \gamma^2}}
\end{equation*}
as follows.  Let $e_{0,2}^* := \lim_{s \rightarrow 2^+} \sum_{\gamma \in \Gamma \setminus \{ 0 \}}%
\overline{\gamma}^2 |\gamma|^{-2s}$.  
Then $\theta(z)$ is given as
$$
	\theta(z) = \Exp{\frac{- e_{0,2}^* z^2 }{2}} \sigma(z).
$$
This function is known to satisfy the transformation formula 
$$
	\theta(z+ \gamma) = \varepsilon(\gamma) \Exp{  \frac{z \overline \gamma}{A} + \frac{\gamma \overline \gamma}{2A}} \theta(z)
$$
for any $\gamma \in \Gamma$, where 
$\varepsilon : \Gamma \rightarrow \{ \pm 1 \}$ is such that
$\varepsilon(\gamma) = -1$ if $\gamma \in 2 \Gamma$ and $\varepsilon(\gamma)=1$ otherwise.
%Hence this function is characterized as the reduced (or normalized) theta function 
%associated to the divisor $[0]$ of the elliptic curve $\bbC/\Gamma$ with 
%$\theta'(0)=1$.  
We define the Kronecker theta function $\Theta(z,w)$ by
$$
	\Theta(z,w) := \frac{\theta(z+w)}{ \theta(z)\theta(w)}.
$$
The above function is known to be a reduced theta function associated to the Poincar\'e bundle
on $\bbC/\Gamma\times\bbC/\Gamma$.

For any $z$, $w \in \bbC$ such that $z$, $w \not\in \Gamma$, we let $K_a(z,w,s) := K^*_a(z,w,s)$, which
we view as a $\mathscr C^\infty$ function for $z$ and $w$.
The relation between this function and the Kronecker theta function is given by
the following theorem due to Kronecker. 
\begin{theorem}[Kronecker]\label{theorem; kronecker}
	$$
		\Theta(z,w) = \exp\left[ \frac{z \overline w}{A}\right] K_1(z,w,1).
	$$
\end{theorem}	
The above theorem was originally proved in terms of Jacobi theta functions by Kronecker
using moduli arguments (See for example \cite{We1}.) 
In \cite{BK1} or \cite{BK2}, we give another proof valid for a fixed lattice $\Gamma\subset\bbC$
using the fact that both sides of the equality are reduced meromorphic theta functions associated to the 
Poincar\'e bundle on $\bbC/\Gamma \times \bbC/\Gamma$,
with the same poles and the same residue at each pole.
%Hence the subtraction of both sides defines a global section of 
%the Poincar\'e bundle, which should be equal to zero since 
%the Hermite form associated to the Poincar\'e bundle is not positive definite. 

%%%%%%%%%%%%%%%%%%%%%%%%%%%%%%%%%%%%%%%%%%%%%%%%%%%%%%%
%
\subsection{Proof of the second limit formula.}
%	
%%%%%%%%%%%%%%%%%%%%%%%%%%%%%%%%%%%%%%%%%%%%%%%%%%%%%%%

In this subsection, we deduce Theorem \ref{thm: KLF} (ii) from Theorem \ref{theorem; kronecker}.
	
\begin{proposition}\label{pro: up to C}
	There exists a constant $C$ such that  
	$$
		\log |\theta(z)|^2-\frac{|z|^2}{A}=-A K^*_0(0,z,1)+C
	$$
	for any $z \notin \Gamma$.
\end{proposition} 

\begin{proof}
	By Theorem \ref{theorem; kronecker} and the fact that 
	$$
		\lim_{z \rightarrow 0} \left[ K_1(z,w,1) -\frac{1}{ z} \right] = K^*_1(0,w,1),
	$$
	we have
	\begin{equation*}
		\lim_{z \rightarrow 0} \left( \Theta(z,w) - \frac{1}{z} \right) 
		= K^*_1(0,w,1) + \frac{\overline w}{A}.
	\end{equation*}
	Direct computation also shows that   
	\begin{equation*}
		\lim_{z \rightarrow 0} \left(\Theta(z,w) - \frac{1}{z} \right)= \frac{\theta'(w)}{\theta(w)}.
	\end{equation*}
	Hence we have
	$$
		K^*_1(0,w,1) + \frac{\overline w}{A} =  \frac{\theta'(w)}{\theta(w)}.
	$$
	In particular, 
	\begin{equation*}\label{equation: leff}
		\frac{\partial}{\partial z} \left( \log \theta(z)  - \frac{z \overline z}{A} \right)= K^*_1(0,z,1).
	\end{equation*}
	Therefore, if we let $\Xi(z)$ be the function 
	$$
		\Xi(z) := \log |\theta(z)|^2 - \frac{|z|^2}{A},
	$$
	then we have
	\begin{align*}
		\frac{\partial}{\partial z} \Xi(z) =  K^*_1(0,z,1), \qquad	\frac{\partial}{\partial \overline z} \Xi(z) 
		=  \overline{K^*_{1}(0,z,1)}. 
	\end{align*}
	On the other hand, one can directly show that 
	\begin{align*}
		A\frac{\partial}{\partial z} K^*_0(0,z,1) =-  {K^*_1(0,z,1)}, \qquad 
		A\frac{\partial}{\partial \overline z} K^*_0(0,z,1)=- {\overline{K^*_1(0,z,1)}}.  
	\end{align*}
	(See for example, Lemma 2.5 and the first formula of p.22 of \cite{BKT}. )
	Hence  $\Xi(z)+AK^*_0(0,z,1)$ must be constant. 
\end{proof}

Our goal is to determine the constant $C$.   We use the following result, which is a type of distribution relation.

\begin{lemma}\label{lem: C}
	We have  
	$$
		\sum_{z_n \not=0} \;K^*_0(0, z_n, 1)=-\frac{2\log n}{A},
	$$
	where the sum is over all $n$-torsion points $z_n$ of $\bbC/\Gamma$ except zero.  
\end{lemma} 

\begin{proof}
	We have 
	$$
		\sum_{z_n \in \frac{1}{n}\Gamma/\Gamma} \pair{\gamma, z_n}=
		\begin{cases}
			n^2 \qquad &(\gamma \in n\Gamma) \\
			0 \qquad &(\gamma \notin n\Gamma). 
		\end{cases}
	$$
	Hence 
	$$
		\frac{1}{n^2}\sum_{z_n} \;K^*_0(0, z_n, s)=\sum_{\gamma \in n\Gamma} 
		\frac{1}{|\gamma|^{2s}}=\frac{1}{n^{2s}}K^*_0(0, 0, s)
	$$
	when the real part of $s$ is sufficiently large, and hence for any $s$ by analytic continuation.  
	In particular, we have
	$$
		\frac{1}{n^2}\sum_{z_n\not=0} \;K^*_0(0, z_n, s)=\left(\frac{1}{n^{2s}}-\frac{1}{n^2}\right)K^*_0(0, 0, s).
	$$
	Since the residue of $K^*_0(0, 0, s)$ at $s=1$ is $1/A$, we have 
	$$
		\frac{1}{n^2}\sum_{z_n\not=0} \;K^*_0(0, z_n, 1)=-\frac{2 \log n}{n^{2}A}
	$$
	as desired.
\end{proof} 
The above lemma shows that the constant $C$ is 
$$
	C=\frac{1}{n^2-1}\left[\sum_{z_n\not=0} \left( \log |\theta(z_n)|^2-\frac{|z_n|^2}{A}\right)-2\log n\right].
$$
We will now calculate this value explicitly in terms of $\Delta$. 

\begin{proposition}\label{pro: C}
	We have 
	$$
		\frac{1}{4} \log|\Delta'|^2=-
		\sum_{z_2\not=0}\left( \log |\theta(z_2)|^2-\frac{|z_2|^2}{A}\right) 
	$$
	where $z_2$ runs through non-trivial $2$-torsion points of $\bbC/\Gamma$ and 
	$$
		\Delta'=(e_1-e_2)^2(e_2-e_3)^2(e_3-e_1)^2
	$$ 
	for $y^2=4x^3-g_2x-g_3=4(x-e_1)(x-e_2)(x-e_3)$.
\end{proposition}

\begin{proof}
	Note that 
	$$
		(x-e_1)(x-e_2)(x-e_3)=\prod_{z_2\not=0}(x-\wp(z_2)).
	$$
	Then if $\Gamma=\bbZ \omega_1+\bbZ \omega_2$,
	we may suppose that $e_1=\wp(\omega_1/2)$, $e_2=\wp(\omega_2/2)$ and 
	$e_3=\wp((\omega_1+\omega_2)/2)$. 
	Since 
	$$
		{\theta(z+w) \theta(z-w)}{\theta(z)^{-2}\theta(w)^{-2}}=\wp(w)-\wp(z),
	$$
	we have  
	$$
		{\theta\left(\frac{\omega_1+\omega_2}{2}\right)\theta\left(\frac{\omega_1-\omega_2}{2}\right)}
		{\theta\left(\frac{\omega_1}{2}\right)^{-2}\theta\left(\frac{\omega_2}{2}\right)^{-2}}
		=e_2-e_1,
	$$
	$$
		{\theta\left(\omega_1+\frac{\omega_2}{2}\right)\theta\left(\frac{\omega_2}{2}\right)}
		{\theta\left(\frac{\omega_1+\omega_2}{2}\right)^{-2}\theta\left(\frac{\omega_1}{2}\right)^{-2}}
		=e_1-e_3,
	$$
	$$
		{\theta\left(\omega_2+\frac{\omega_1}{2}\right)\theta\left(\frac{\omega_1}{2}\right)}
		{\theta\left(\frac{\omega_1+\omega_2}{2}\right)^{-2}\theta\left(\frac{\omega_2}{2}\right)^{-2}}
		=e_2-e_3.
	$$
	Hence using the transformation formula of $\theta(z)$, the value $\Delta'$ is  
	$$
		\exp\left[\frac{\omega_1\overline{\omega_1}+\omega_2\overline{\omega_2}+\overline{\omega_1}\omega_2}{A}\right] 
		{\theta\left(\frac{\omega_1}{2}\right)^{-4}
		\theta\left(\frac{\omega_2}{2}\right)^{-4}\theta\left(\frac{\omega_1+\omega_2}{2}\right)^{-4}}.
	$$
	Multiplying it and its complex conjugation and taking the logarithm, 
	we obtain the formula. Note that since we take the 
	logarithm of  {\it positive real} numbers, the values do not depend on the choice of the branch of the logarithm.
\end{proof}

\begin{proof}[Proof of Theorem \ref{thm: KLF} (ii)]
	Since the Ramanujan $\Delta$ is given by $\Delta=2^4 \Delta'$, we have by Lemma \ref{lem: C} ad Proposition \ref{pro: C}
	$$
		C=\frac{1}{3}\left(-\frac{1}{4} \log|\Delta'|^2-2\log 2\right)=-\frac{1}{12}\log|\Delta|^2.
	$$
	Our assertion now follows from Proposition \ref{pro: up to C}.
\end{proof}

%%%%%%%%%%%%%%%%%%%%%%%%%%%%%%%%%%%%%%%%%%%%%%%%%%%%%%%
%
\subsection{Proof of  the first limit formula.} 
%
%%%%%%%%%%%%%%%%%%%%%%%%%%%%%%%%%%%%%%%%%%%%%%%%%%%%%%%

We now prove Theorem \ref{thm: KLF} (i) using the second limit formula.

\begin{proof}[Proof of Theorem \ref{thm: KLF} (i)]
	From (\ref{eq: integral expression}), we have 
	\begin{align*}
		A^{s-1} \Gamma(s) & \left(A K^*_0(0, 0, s)-\frac{1}{s-1}\right)\\
		&= I_0(0, 0, s)
	 	- \frac{1}{s} + I_{0} (0, 0, 1-s)  - \frac{A^{s-1} \Gamma(s)-1}{s-1}.
	\end{align*}
	Therefore, we have
	\begin{equation}\label{eq: integral expression 1}
		\lim_{s \rightarrow 1}\left(A K^*_0(0, 0, s)-\frac{1}{s-1}\right) = I_0(0, 0, 1)
	 	-1 + I_{0} (0, 0, 0)
		-\log A+c,
	\end{equation}
	where $c$ is the Euler constant as before and we used the fact $\Gamma'(1)=-c$. 
	On the other hand, we have 
	\begin{equation*}
		 AK^*_0(0, z, 1) = I_0(0, z, 1) - 1 + I_{0} (z, 0, 0).
	\end{equation*}
	We let $$I_{0}^* (z, 0, s)= I_{0} (z, 0, s)-\int_{1}^\infty \exp(-t|z|^2/A) t^{s-1} dt.$$
	Then 
	$\displaystyle\lim_{z \rightarrow 0}  I_0(0, z, 1)= I_0(0, 0, 1)$ and
	$\displaystyle\lim_{z \rightarrow 0}  I_0^*(z, 0, 0)= I_0(0, 0, 0)$.
	We have 
	\begin{align*}
		\Gamma(s)-\frac{1}{s}&=
		\int_{|z|^2/A}^\infty e^{-t} t^{s-1} dt+\int_{0}^{|z|^2/A} e^{-t} t^{s-1} dt-\frac{1}{s}\\
		&=
		\int_{|z|^2/A}^\infty e^{-t} t^{s-1} dt+\int_{0}^{|z|^2/A} (e^{-t}-1) t^{s-1} dt+
		\frac{1}{s}\left[ \left(\frac{|z|^2}{A}\right)^s-1\right].
	\end{align*}
	Taking $s \rightarrow 0$, we have 
	$$
		-c=\int_{|z|^2/A}^\infty e^{-t} t^{-1} dt+\int_{0}^{|z|^2/A} (e^{-t}-1) t^{-1} dt
		+\log \left(\frac{|z|^2}{A}\right).
	$$
	Hence 
	\begin{align*}
	 	&AK^*_0(0, z, 1) = I_0(0, z, 1) - 1 + I^*_{0} (z, 0, 0)+\int_{1}^\infty \exp(-t|z|^2/A) t^{-1} dt \\
		&= I_0(0, z, 1) - 1 + I^*_{0} (z, 0, 0)
		-c-\int_{0}^{|z|^2/A} (e^{-t}-1) t^{-1} dt-\log \left(\frac{|z|^2}{A}\right).
	\end{align*}
	Therefore 
	\begin{align*}
		 \lim_{z \rightarrow 0} 
		\left(AK^*_0(0, z, 1)+\log |z|^2 \right)= I_0(0, 0, 1) - 1 + I_{0} (0, 0, 0)
		-c+\log A.
	\end{align*}
	Finally, combining this with (\ref{eq: integral expression 1}) and the second limit formula, we have 
	\begin{align*}
		\lim_{s \rightarrow 1}\left(A K^*_0(0, 0, s)-\frac{1}{s-1}\right)&= 
		 \lim_{z \rightarrow 0} 
		\left(AK^*_0(0, z, 1)+\log |z|^2 \right) -2\log A+2c \\
		=&-\frac{1}{12}\log|\Delta|^2-2\log A+2c. 
	\end{align*}
	This proves our assertion.
\end{proof}

%%%%%%%%%%%%%%%%%%%%%%%%%%%%%%%%%%%%%%%%%%%%%%%%%%%
%
%
%
%
\section{$p$-adic Kronecker  limit formulas}
%
%
%
%%%%%%%%%%%%%%%%%%%%%%%%%%%%%%%%%%%%%%%%%%%%%%%%%%%

%%%%%%%%%%%%%%%%%%%%%%%%%%%%%%%%%%%%%%%%%%%%%%%%%%%
%
\subsection{The $p$-adic Eisenstein-Kronecker functions}
%
%%%%%%%%%%%%%%%%%%%%%%%%%%%%%%%%%%%%%%%%%%%%%%%%%%%

Assume now the conditions of the second half of the introduction.
In \cite{BFK}, we defined a $p$-adic analogue of the  Kronecker double  series as a Coleman function on a
CM elliptic curve. 
In order to prove the $p$-adic limit formulas, we 
define in this subsection a $p$-adic analogue of the function $\log |\theta(z)|^2-|z|^2/A$, 
which turns out to be a Coleman function.   We then prove the distribution relation, which
will be used to characterize this function.

Let $p$ be a prime $\geq 5$.  In what follows, fix an embedding of $\overline{\bbQ}$ into $\bbC_p$. 
We fix a branch of the logarithm, which is a homomorphism $\log_p: \bbC_p^\times \rightarrow \bbC_p$.  We extend this homomorphism 
to $\bbC_p[[t]]$ by using the decomposition $\bbC_p[[t]]=\bbC_p \times (1+t\bbC_p[[t]])$ and defining
 $\log_p (1-tf(t))=-\sum {t^nf^n(t)}/n$ for any $f(t) \in \bbC_p[[t]]$. 
Let $E$ be a CM elliptic curve as in the introduction and let $\Gamma$ be the period lattice of $E\otimes\bbC$. 
For $z_0 \in \Gamma \otimes \bbQ$, we let 
  $$
    \theta_{z_0}(z):=\theta(z+z_0) \exp\left(-\frac{z\overline{z_0}}{A}-\frac{z_0\overline{z_0}}{2A}\right).
  $$
Then by \cite{BK1}, the Taylor  series of $\theta_{z_0}(z)$ at $z=0$ has algebraic coefficients.
If we consider the formal composition 
$$
    \widehat{\theta}_{z_0}(t):=\theta_{z_0}(z)|_{z=\lambda(t)}
$$
of this series with $\lambda(t)$, where $\lambda(t)$ is the formal logarithm of the formal group of $E$, then we may regard this 
power series as an element in  $\bbC_p[[t]]$. 
Considering its derivatives, we may prove that 
$\log_p \widehat{\theta}_{z_0}(t)$ is a rigid analytic function on 
the open unit disc over $\bbC_p$, namely, it is convergent if $|t|<1$. 

We use the same notations as in \cite{BFK}.   In particular, we fix a prime $\frp$ in $\cO_\bsK$ over $p \geq 5$,
and we let $\pi:= \psi_{E/\bsK}(\frp)$, where $\psi_{E/\bsK}$ is the Gr\"ossen character of $K$ associated to the elliptic curve $E$.
Then $\pi$ is a generator of the ideal $\frp$.

\begin{definition}\label{def: log-p-theta}
	We let $\log _p \theta$ be the function in $A_{\log}(U)$ defined by 
	$$
    		\log _p \theta |_{]z_0[}:=\log_p \widehat{\theta}_{z_0}(t) \in A_{\log}(]z_0[)
	$$
	on each residue disc $]z_0[$. 
\end{definition}

Now we investigate basic properties of $\log_p {\theta}$. 

\begin{proposition} 
	For $z_0 \in \Gamma \otimes \bbQ$ and $z_\alpha$ such that 
	$\alpha z_\alpha \in \Gamma$ for $\pi$-power morphism $\alpha \in \mathrm{End}_{\overline{\bbQ}}(E)$, we have 
	$$
  		\log_p \widehat{\theta}_{z_0}(t \oplus t_\alpha)=\log_p\widehat{\theta}_{z_0+z_\alpha}(t)
 	$$
	where $z_\alpha \in \bbC$ is a lift of a torsion point in $\bbC/\Gamma$ 
	corresponding to $t_\alpha \in E(\bbC_p)_{\mathrm{tor}}$, and the right hand side is independent 
	of the choice of the lift. 
\end{proposition}

\begin{proof}
          Let $\alpha$ and $\beta$ be elements of $\cO_K$ such that 
          $2\alpha | \beta$ and $\beta z_0 \in \Gamma$. 
           Then  $f_{\beta}(z):=\theta(z)^{N\beta}/\theta( \beta z)$  is a rational function on $E$ over $\overline{\bbQ}$.  
           We have 
 	\begin{equation}\label{equation: logtheta}
    		\theta_{z_0+z_\alpha}(z)^{N\beta}= \pm \theta( \beta z) 
 		\tau_{z_0+z_\alpha}^*f_{\beta}(z). 
 	\end{equation}
	Similarly,  we have 
  	$$
        		\widehat{\theta}_{z_0}(z)^{N\beta}
                  = \pm \theta( \beta z) \tau_{z_0}^*f_{\beta}(z). 
  	$$
	Since  $f_\beta$ is a rational function, we also have 
  	$$
      		\tau_{z_0+z_\alpha}^*f_{\beta}(t)=  \tau_{z_0}^*f_{\beta}(t\oplus t_\alpha).
  	$$
	Hence  we have 
   	\begin{equation}\label{equation: logtheta2}
    		\widehat{\theta}_{z_0}(t \oplus t_\alpha)^{N\beta}
                  = \pm \theta( [\beta] t) \tau_{z_0+z_\alpha}^*f_{\beta}(t). 
  	\end{equation}
	Our assertion now follows from \eqref{equation: logtheta} and \eqref{equation: logtheta2}.
\end{proof}

\begin{corollary}\label{corollary: interpolation}
	Let $t_\alpha$ be a $\pi$-power torsion point, and we assume that $z_0 \not=0$ or $t_\alpha\not=0$. 
	Then we have 
	$$
  		\log_p \widehat{\theta}_{z_0}(t_\alpha)=
		\log_p \left(\widehat{\theta}(z_0+z_\alpha) \exp\left[-\frac{(z_0+z_\alpha)\overline{(z_0+z_\alpha)}}{2A}\right] \right).
	$$
\end{corollary}

Roughly speaking,  $\log_p \theta(z)$ is a $p$-adic function which interpolates the special values $\log \theta(z)-z \overline{z}/2A$ at torsion points. 
We thus regard  $\log_p \theta(z)$ as a $p$-adic analogue of the function $\log |\theta(z)|^2-|z|^2/A$.

%%%%%%%%%%%%%%%%%%%%%%%%%%%%%%%%%%%%%%%%%%%%%%%%%%%
%
\subsection{The $p$-adic second limit formula}
%
%%%%%%%%%%%%%%%%%%%%%%%%%%%%%%%%%%%%%%%%%%%%%%%%%%%

In this subsection, we prove Theorem \ref{thm: pKLF}, which is a $p$-adic analogue of the Kronecker
second limit formula.

\begin{proposition}\label{proposition: theta distribution}
         Let $z_0$ be a $\mathfrak{f}$-torsion point of $\bbC/\Gamma$. Then for $\alpha \in \cO_K$ 
	we have 
	$$
      	\theta_{\alpha z_0}( \alpha z)^{24 N(\alpha\mathfrak{f})}= \Delta^{2N(\alpha\mathfrak{f})(N\alpha-1)} 
                        \prod_{z_\alpha \in E[\alpha]} \theta_{z_0+z_\alpha}(z)^{24 N(\alpha\mathfrak{f})}
	$$
where $z_\alpha$ is a lift of a $\alpha$-torsion point of $E$ and 
the right hand side is independent of the choice of the lifts $z_0$ and $z_\alpha$ on $\bbC$. 
\end{proposition}
\begin{proof}
         Since for $\gamma \in \Gamma$  we have 
         $\theta_{z_0+\gamma}(z)=\pm \pair{z_0/2, \gamma} \theta_{z_0}(z)$, 
         the function $ \theta_{z_0}(z)^{2N\mathfrak{f}}$ is independent of the lift $z_0$ if 
         $(N\mathfrak{f})z_0 \in \Gamma$. The independence of the lifts of $z_0$ and $z_\alpha$ follows from this fact. 
		The logarithmic  derivatives of both sides coincide (see for example Proposition 2.15 of \cite{BFK}). 
	Hence for each $\alpha$,  there exists a constant $c_\alpha$ such that 
	 $$
     		\theta_{\alpha z_0}( \alpha z)^{2N(\alpha\mathfrak{f})}= c_\alpha 
                        \prod_{z_\alpha \in E[\alpha]} \theta_{z_0+z_\alpha}(z)^{2N(\alpha\mathfrak{f})}. 
	  $$
	By the definition of $\theta_{z_0}(z)$, we see that $c_{\alpha}$ is independent of  $z_0$. 
	Hence we may assume that $z_0=0$ and $N\mathfrak{f}=1$. 
	Then we have 
	\begin{align*}
		 \prod_{z_{\alpha \beta} \in E[\alpha \beta]} \theta_{z_{\alpha \beta}}(z)^{2N(\alpha\beta)}
		= \prod_{z_{\alpha \beta} \in E[\alpha \beta]/E[\alpha]} \prod_{z_{\alpha } \in E[\alpha ]} 
		  \theta_{z_{\alpha \beta}+z_{\alpha}}(z)^{2N(\alpha\beta)}\\
		= \prod_{z_{\alpha \beta} \in E[\alpha \beta]/E[\alpha]}
		c_\alpha^{N\beta}  \theta_{ \alpha z_{\alpha \beta}}(\alpha z)^{2N(\alpha\beta)}\\
		= c_\alpha^{N\beta^2} c_\beta^{2N\alpha} \theta(\beta \alpha z)^{2N(\alpha\beta)}.
	\end{align*} 
	Hence we have 
 	$
   		c_\alpha^{N\beta^2} c_\beta^{2N\alpha}=c_{\alpha \beta}=c_\beta^{N\alpha^2} c_\alpha^{2N\beta}
 	$
	or equivalently,  
	$$
   		c_\alpha^{N\beta(N\beta-1)}=c_\beta^{N\alpha(N\alpha-1)}. 
 	$$
	In particular, $c_\alpha^{12}=c_2^{N\alpha(N\alpha-1)}$. 
	On the other hand, we consider the constant term of 
	$$
   		\frac{\theta(2z)^8}{\theta(z)^8}=c_2   \prod_{z_2 \in E[2]-\{0\}} \theta_{z_2}(z)^{8}. 
	$$
	As in the proof of Proposition \ref{pro: C}, 
	we have 
	$$
  		\prod_{z_2 \in E[2]-\{0\}} \theta_{z_2}(0)^{8}=\Delta'^{-2}.
	$$
	Hence $c_2=2^8 \Delta'^2=\Delta^2$. Our assertion now follows from these facts. 
\end{proof}

\begin{corollary}
	The function $\Xi(z):= - \log_p \theta (z) -\frac{1}{12}\log_p \Delta$ satisfies the distribution relation  
	$$
		\Xi(\alpha z)=\sum_{z_\alpha \in E[\alpha]} \Xi(z+z_\alpha). 
	$$
\end{corollary}

\begin{proof}
	By Proposition \ref{proposition: theta distribution}, we have 
	$$
		\log_p \widehat{\theta}_{\alpha z_0}([\alpha]t)=\frac{N\alpha-1}{12} \log_p \Delta 
		+\sum_{t_\alpha \in E[\alpha]} \log_p \widehat{\theta}_{z_0}(t \oplus t_\alpha). 
	$$
	Our assertion follows from this formula. 
\end{proof}

We now prove the $p$-adic second limit formula.

\begin{proof}[Proof of Theorem \ref{thm: pKLF}]
	Since the derivatives of $\Xi(z)$ and 
	$$
		K_0^\col(0,z,1):=E^{\col}_{1,1}(z)
	$$ are equal, their difference 
	$c(z)=\Xi(z)-E^{\col}_{1,1}(z)$ is a constant on  the residue disc $]z_0[$. 
	By (i) and the definition of $E^{\col}_{1,1}(z)$, the locally constant function $c(z)$ satisfies the distribution relation. 
	For any torsion point $z_0$ of order $\mathfrak{f}$, we take $N$ such that $\pi^N \equiv 1 \mod \mathfrak{f}$. 
	Then $[\pi^N]^*(]z_0[)=]z_0[$ and 
	$$
		[\pi^N]^*c(z) |_{]z_0[}=\sum_{w \in E[\pi^N]} c(z+w) |_{]z_0[}. 
	$$
	Since $c(z) |_{]z_0[}$ is constant, the above relation shows $c(z) |_{]z_0[}=0$.  
\end{proof}

The above result shows in particular that  $\Xi(z)=- \log_p \theta (z) -\frac{1}{12}\log_p \Delta$ is in fact a Coleman function.

%%%%%%%%%%%%%%%%%%%%%%%%%%%%%%%%%%%%%%%%%%%%%%%%%%%
%
\subsection{$p$-adic Eisenstein-Kronecker-Lerch series}
%
%%%%%%%%%%%%%%%%%%%%%%%%%%%%%%%%%%%%%%%%%%%%%%%%%%%

In this subsection, we introduce the $p$-adic Eisenstein-Kronecker-Lerch series.  Then
we consider a $p$-adic counterpart to the arguments concerning the Kronecker first limit formula 
in the classical case and prove Proposition \ref{pro: pKLF}.

Let $p \geq 5$ be a prime of good \textit{ordinary} reduction for $E$, and we fix a prime $\frp$ of $\cO_\bsK$ over $p$.  
We defined in \cite{BK1} \S 3.1 a $p$-adic measure $\mu:=\mu_{0, 0}$ on $\bbZ_p \times \bbZ_p$ interpolating
the Eisenstein-Kronecker numbers, or more precisely,
the special values of Eisenstein-Kronecker-Lerch series $K^*_{a+b}(0,0,b; \Gamma)/A(\Gamma)^a$ for $a$, $b \geq 0$,
where $\Gamma$ is the period lattice of $E$.
We define the $p$-adic Eisenstein-Kronecker-Lerch function as in the introduction as follows.

\begin{definition}  For any integer $a \in \bbZ$, we define the \textit{$p$-adic Eisenstein-Kronecker-Lerch function} by
	$$
		K^{(p)}_a(0,0,s) := \int_{\bbZ_p^\times \times \bbZ_p^\times} \pair{x}^{s-1}\pair{y}^{a-s} \omega(y)^{a-1} d \mu(x,y).
	$$
\end{definition}

%\begin{theorem}
%We have 
%$$\int_{\bbZ_p^\times \times \bbZ_p^\times} \langle x \rangle^{s-1} \langle y\rangle^{-s} d\mu(x,y)=$$
%where $s$ in the limit is taken  over $s \in (p-1)\bbZ$ which is $p$-adically converges to $1$. 
%$$\int_{\bbZ_p^\times \times \bbZ_p^\times}y^{-1} d\mu(x,y)=\left(1-\frac{1}{p}\right) \log_p \overline{\pi}.$$
%\end{theorem}

The $p$-adic Eisenstein-Kronecker-Lerch function is analytic in $s \in \bbZ_p$. 
The reason we view this function as a $p$-adic analogue of Eisenstein-Kronecker-Lerch series is the
following interpolation property.

\begin{proposition}
	For any integer $a$, $b$ such that $a \geq b > 0$ and $b \equiv 1 \pmod{p-1}$, we have
	\begin{equation}\label{eq: interpolation}%\begin{multline*}
		\frac{K^{(p)}_a(0,0,b)}{\Omega_p^{a-1}} 
		=(-1)^{a-1}(b-1)! \left(1-\frac{\pi^{a}}{p^{a-b+1}}\right)
		\left(1-\frac{\pi^{a}}{p^b} \right)\frac{K^*_{a}(0,0,b)}{A(\Gamma)^{a-b}},
	\end{equation}%\end{multline*}
	where $\Omega_p$ is a $p$-adic period of the formal group of $E$.
\end{proposition}

\begin{proof}
	This follows from the interpolation property of the measure $\mu:=\mu_{0,0}$ given in \cite{BK1} Proposition 3.5.
\end{proof}

We now give the proof of Proposition \ref{pro: pKLF}.  %For simplicity, we let $E_{1,1}(z; \Gamma) := K^\col_0(0,z,1; \Gamma)$.

\begin{proof}[Proof of Proposition \ref{pro: pKLF}]
	We consider the function 
	$$
		f(t):=\Omega_p \int_{\bbZ_p^\times \times \bbZ_p^\times}y^{-1} \exp(y\Omega_p^{-1} \lambda(t)) d\mu(x,y)
	$$
	on the $p$-adic residue disc $]0[$ around $0$.
	Then 
	\begin{multline*}
 		\lambda'(t)^{-1}\frac{d}{dt}f(t)=\int_{\bbZ_p^\times \times \bbZ_p^\times} \exp(y\Omega_p^{-1} \lambda(t)) d\mu(x,y)\\
		= F_1(t; \Gamma)-{\overline{\pi}}^{-1}F_1([\pi]t; \Gamma)-F_1(t; \ol\frp\Gamma)
		 +{\overline{\pi}}^{-1} F_1([\pi]t; \ol\frp\Gamma).
	\end{multline*}
	Hence for $E_{1,1}^{(p)}(z; \Gamma) := E_{1,1}^\col(z; \Gamma) - p^{-1} E_{1,1}^\col(\pi z; \Gamma)$, the function
	$$
		f(t)-E_{1,1}^{(p)}(z; \Gamma)+E_{1,1}^{(p)}(z; \ol\frp\Gamma)
	$$
	is a constant on the residue disc $]0[$. 
	Furthermore, both functions satisfy the distribution relation $\sum_{t_{\pi} \in E[\pi]} f(t\oplus t_{\pi})=0$ and 
	$\sum_{t_{\pi} \in E[\pi]} E_{1,1}^{(p)}(t\oplus t_{\pi})=0$. Hence we must have 
	$f(t)=E_{1,1}^{(p)}(z; \Gamma)-E_{1,1}^{(p)}(z; \ol\frp\Gamma)$ 
	on ${]0[}$. On the other hand, the $p$-adic  second limit formula shows that 
	$$
		E_{1,1}^{(p)}(z; \Gamma)=\log_p \theta(z; \Gamma)
		-\frac{1}{p}\log_p \theta(\pi z; \Gamma)+\frac{1}{12}\left(1-\frac{1}{p}\right) \log_p \Delta(\Gamma), 
	$$
	hence we have
	$$
		f(0) = \bigl.E_{1,1}^{(p)}(z; \Gamma)-E_{1,1}^{(p)}( z; \overline{\pi}\Gamma)\;\bigr|_{z=0}=
		\left(1-\frac{1}{p}\right) \log_p \overline{\pi}.
	$$
	Our assertion now follows from the fact that
	$
		f(0)=\Omega_p K^{(p)}_0(0,0,1).
	$
\end{proof}

\begin{remark}\label{rem: not pKLF}
	In the interpolation formula of \eqref{eq: interpolation}, if we let $a=0$ and $b=1$,  then
	the interpolation factor of the the right hand side vanishes. Hence the value 
	$$
		\Omega_p K^{(p)}_0(0,0,1) = \int_{\bbZ_p^\times \times \bbZ_p^\times}  y^{-1} d\mu(x,y)
	$$
	is in some sense not the constant term but 
	the residue at $s=1$ of the $p$-adic analogue of $\sum_{\gamma \in \Gamma}^*1/|\gamma|^{2s}$. 
	Because of this fact, the formula of Proposition \ref{pro: pKLF} is not a perfect $p$-adic analogue of the classical Kronecker first limit formula. 
\end{remark}

%%%%%%%%%%%%%%%%%%%%%%%%%%%%%%%%%%%%%%%%%%%%%%%%%%%%%%%%%%%%%%%%%%


\begin{thebibliography}{}
  	\bibitem[BK1]{BK1}{K. Bannai and S. Kobayashi}, Algebraic theta functions and 
	$p$-adic interpolation of Eisenstein-Kronecker numbers, preprint,  arXiv:math/0610163v4 [math.NT], 
	11 December 2007. %, 52 pages.
	\bibitem[BK2]{BK2}{K. Bannai and S. Kobayashi}, Algebraic theta functions and
	Eisenstein-Kronecker numbers, RIMS K\^oky\^uroku Bessatsu B4: \textit{Proceedings of the Symposium on Algebraic Number theory and Related Topics},  
	eds. K. Hashimoto, Y. Nakajima and H. Tsunogai, December (2007), 63--78.
	\bibitem[BKT]{BKT}{K. Bannai, S. Kobayashi and T. Tsuji}, On the real Hodge and $p$-adic realizations of the elliptic polylogarithm for CM elliptic curves,
	Preprint 2007, arXiv:0711.1701v1 [math.NT].
	\bibitem[BFK]{BFK}{K. Bannai, H. Furusho and S. Kobayashi}, in preparation.
	\bibitem[dS]{dS}{E. de Shalit}, \textit{Iwasawa theory of elliptic curves with complex multiplication},
	Academic Press (1987). 
	\bibitem[Ka]{Ka}N. Katz, $p$-adic interpolation of real analytic Eisenstein
	series, Ann.\ of Math.\ \textbf{104} (1976), 459-571.
	\bibitem[We1]{We1}{A. Weil}, {\it Elliptic Functions according to Eisenstein and Kronecker},
         Springer-Verlag, 1976. 
\end{thebibliography}
\end{document}